\documentclass{rspublic}
\usepackage{amssymb, amsmath, amsthm, mathrsfs}

\newtheorem{statement}[theorem]{Statement}

\def\Re{{\mathbb R}}

\def\d{{\rm d}}

\def\d{{\rm d}}

\def\e{{\rm{e}}}

\def\R{\mathbb R}
\def\lap{\Delta}
\def\grad{\nabla}
\def\div{\nabla\cdot}

\def\Re{\R}

\def\qqand{\qquad\mbox{and}\qquad}

\def\be#1{\begin{equation}\label{#1}}
\def\ee{\end{equation}}
\def\bea{\begin{eqnarray*}}
\def\eea{\end{eqnarray*}}

\def\({\left(}
\def\){\right)}

\title[Regularity via numerics]
{Numerical verification of regularity in the three-dimensional Navier-Stokes
equations}
\author[J.C. Robinson \& W. Sadowski]
{James C. Robinson$^1$ \& Witold Sadowski$^{2,3}$}
\affiliation{$^1$Mathematics Institute, University of Warwick,
Coventry CV4 7AL.\\ UK.\\ $^2$ Faculty of Mathematics, Informatics
and Mechanics,\\
       Warsaw University,
       Banacha 2 02-097, Warszawa, Poland. \\ $^3$ Institute of Mathematics of the Polish Academy of
       Sciences,\\
       ul. \'Sniadeckich 8 00-956, Warszawa, Poland.\\}

\begin{document}

\maketitle

\begin{abstract}{Navier-Stokes equations, rigorous computation, Galerkin method}
Current theoretical results for the three-dimensional Navier--Stokes equations
 only guarantee that solutions remain regular for all time when the initial
 enstrophy ($\|Du_0\|^2:=\int|{\rm curl}\,u_0|^2$) is sufficiently small,
 $\|Du_0\|^2\le\chi_0$. In fact, this smallness condition is such that the enstrophy
 is always non-increasing. In this paper we provide a numerical procedure that will
 verify regularity of solutions for any bounded set of initial conditions,
 $\|Du_0\|^2\le\chi_1$. Under the assumption that the equations are in fact regular we show that
 this procedure can be guaranteed to terminate after a finite time.
\end{abstract}

\section{Introduction}

The global regularity of solutions of the three-dimensional
Navier-Stokes equations is one of the most intriguing open
problems in mathematics. Despite the efforts of many eminent
mathematicians over the last 150 years, the question whether the most
fundamental model in hydrodynamics is or is not self consistent
remains unanswered. Even when the flow is not
disturbed by any external forces and the problem of boundaries is avoided by considering periodic boundary conditions (or the whole of $\Re^3$) the problem is unresolved.

More precisely, consider the Navier-Stokes equations for the flow of an
incompressible fluid
\begin{equation}\label{ns}
\frac{\partial u}{\partial t}-\nu  \lap u+(u\cdot\grad)u+\grad
p=0\qquad\div u=0\qquad u(0)=u_0.
\end{equation}
with periodic boundary conditions on $\Omega=[0,L]^3$, where $\nu>0$ is the kinematic
viscosity. It is known that weak solutions of these
equations exist that are defined for all $t\ge0$. However, we do not know whether these solutions are
unique and regular. On the other hand, we know that if the initial enstrophy is finite then the solutions are
regular for $t$ sufficiently small ($t<T_1$) and sufficiently
large ($t>T_2$), where $T_1$ and $T_2$ depend respectively on the
enstrophy and kinetic energy of initial conditions:
\begin{equation}\label{Tbound} T_1 \sim \|Du_0\|^{-4}\qquad{\rm and}
\qquad T_2 \sim \|u_0\|^2,
\end{equation}
where
$$
\|u_0\|^2=\int_\Omega|u(x)|^2\,\d x\qquad\mbox{and}\qquad\|Du_0\|^2=\int_\Omega|{\rm curl}\,u(x)|^2\,\d x.
$$
It is an open
problem whether or not weak solutions develop any singularity at
some time $T_1<t<T_2$.

In the light of (\ref{Tbound}) it is natural to conjecture that some
restrictions on enstrophy or energy of initial conditions might
provide us with regularity of weak solutions. However, nothing is
known about regularity of weak solutions only under the assumption
that energy of the flow is small. In other words, there is no
theoretical evidence that there exists an $\varepsilon
>0$ such that all initial conditions with kinetic energy
$$
\|u_0\|^2< \varepsilon
$$
do not blow up. On the
other hand, under the stronger assumption that the enstrophy of the
initial condition $u_0$ is sufficiently small,
$$
\|Du_0\|^2< \varepsilon _1
$$
one can prove that $u_0$ gives rise to a strong global solution. The value
of $\varepsilon _1$ for which the argument holds is very small (see Section \ref{setting}). This
is due to the fact that the proof of regularity of solutions
arising from small initial data (sketched in the next
section) actually rules out all initial conditions for which the
enstrophy is increasing at the beginning of the flow. Thus it is
desirable to extend this theoretical result to a larger class of
initial conditions that allow for more complicated and physically
realistic behaviour of the flow.

The main aim of this paper is to
present a numerical method which may provide us with such an extension. We show that, under the assumption that the equations are well-posed, it is possible to verify numerically the regularity of solutions for all initial conditions with $\|Du_0\|\le R$ for any choice of $R>0$. In Section \ref{setting} we set up the problem more formally, and recall some classical results on the regularity of solutions. Section \ref{singleu0} reproduces results from Dashti \& Robinson (2006) that can be used to verify that a single initial condition $u_0\in V$ gives rise to a strong solution of some bounded interval $[0,T]$. Section \ref{finitetime} shows that for initial conditions in any bounded ball in $V$ it is sufficient to guarantee that solutions are regular on some finite time interval $[0,T]$, since all solutions are regular after some finite time. Section \ref{H2ball} is the core of the paper, and shows that one can verify that all initial conditions in some bounded ball in $H^2$ give rise to regular solutions -- this requires no additional assumptions. The final section shows his this can be used to verify the same property for all initial conditions in a bounded ball in $H^1$, although we note here that this requires the assumption of regularity for all initial conditions in $H^1$.

\section{Setting of the problem}\label{setting}

As enstrophy and kinetic energy play key role in the problem of
existence and regularity of solutions it is worth recalling two
basic facts concerning them. First of all, it follows from the energy inequality
\begin{equation}
\frac{\d}{\d t} \|u\|^2 + \nu \|Du\|^2 \le 0
\end{equation}
that the kinetic energy of the flow decreases
to zero as time goes to infinity (the flow is dissipative).
Moreover, we have
\begin{equation}\label{boundonDu}
\int_0^T \|Du(s)\|^2\,\d s \leq
\nu^{-1}\|u(0)\|^2
\end{equation}
which shows that mean value of $\|Du\|^2$ over the interval $[0,T]$ is bounded
in terms of the initial kinetic energy and decreasing at least as fast as $T^{-1}$.

Throughout the paper we use standard notation for the Lebesgue and
Sobolev spaces $L^2$, $H^1$ and $H^2$. The set of all smooth,
periodic and divergent free functions on $\Omega$ with $\int_\Omega u(x)\,\d x=0$ we denote by
$\mathcal{V}$. The closure of $\mathcal{V}$ in the $L^2$ and $H^1$
norms we denote by $H$ and $V$, respectively. The eigenfunctions
of the Stokes operator $A$ defined by
$$
A(u)=-P\Delta (u),
$$
where $P$ is the
Leray projector (the orthogonal projection of the space $L^2$ onto $H$), we denote by $w_k$
and their corresponding eigenvalues by $\lambda_k$. It follows
from the theory of compact self-adjoint operators in Hilbert
spaces that we can order the eigenvalues so that
$$
0<\lambda_1 \le \lambda_2 \le \lambda_3 \le ...
\rightarrow \infty.
$$
Since we consider periodic boundary
conditions we can in fact specify the values of $\lambda _k$
exactly, and in particular we have
$$
\lambda_1=\frac{4\pi^2}{L^2}.
$$
In what follows we will concentrate on a case $L=2\pi$ for which
$\lambda_1=1$. This is convenient in the light of the
following Poincar\'e inequalities
$$
\|u\|^2 \le \lambda_1
^{-1} \|Du\|^2,\qquad \|Du\|^2 \le\lambda_1 ^{-1}
\|Au\|^2
$$
which we use frequently in what follows.

As mentioned above, solutions to (\ref{ns}) with sufficiently small
initial data are regular for all time. With zero forcing the
condition guaranteeing regularity is given by following theorem.

\begin{theorem}\label{c}
If the initial condition $u_0$ for problem
(\ref{ns}) satisfies
\begin{equation}\label{small}
\|Du_0\|^2\leq c^{-1/2} \nu ^2\lambda_1^{1/2},
\end{equation}
where $c$ is some absolute constant\footnote{For later use it is
worthwhile to note that $c=27k^4/16$, where $k$ is the constant in the inequality
$$
\left|\int_\Omega [(u\cdot\grad)v]\cdot\Delta w\,\d x\right|\le
k\|Du\|\|Dv\|^{1/2}\|Av\|^{1/2}\|Aw\|.
$$
In the cubic geometry here we can take $k=9\times 2^{15/4}$, see
Dashti \& Robinson (2006).}, then $u_0$ gives rise to a global
strong solution.
\end{theorem}

The proof is based on the inequality
\begin{equation}
\label{enstrophy}
\frac{\d}{\d t}\|Du\|^2 \leq \frac{c}{\nu ^3} \|Du\|^6 -\nu \lambda_1
\|Du\|^2
\end{equation}
which may be obtained from the enstrophy inequality with some help
from the Sobolev embedding theorem. What is crucial is that for
initial conditions satisfying (\ref{small}) the right hand side of
(\ref{enstrophy}) is non-positive. In consequence $\|Du\|^2$ will
never exceed the value $c^{-1/2} \nu ^2\lambda_1^{1/2}$. (For
details of the proof see, for example, Theorem * in Constantin \&
Foias (1988).)

Below we treat the question of regularity for
periodic boundary conditions on $[0,2\pi]^3$ and for viscosity
$\nu=1$. We now show, using a simple scaling argument, that a
proof of regularity for this particular case would imply
regularity for all $[0,L]^3$ and all $\nu>0$. Indeed, starting with
the equation
$$
\frac{\partial u}{\partial t}-\nu\Delta
u+(u\cdot\grad)u+\grad p=0\qquad\nabla\cdot u=0,\qquad
x\in[0,L]^3,
$$
we set
$$
\tilde u=Lu/2\pi\nu,\quad\tilde x=2\pi
x/L,\quad\tilde t=4 \pi ^2\nu t/L^2,\quad{\rm and}\quad\tilde
p=L^2p/4 \pi^2\nu^2,
$$
and then obtain
\begin{equation}\label{ournse}
\frac{\partial\tilde
u}{\partial \tilde t}-\Delta_{\tilde x}\tilde u+(\tilde
u\cdot\grad_{\tilde x})\tilde u+\grad_{\tilde x}\tilde p=0
\end{equation}
for $\tilde x\in[0,2 \pi]^3$.

For this non-dimensionalized form of the Navier-Stokes equations we
can get an explicit value $R_V$ for the radius  of the ball in $V$
in which standard results guarantee that every initial condition
gives rise to a global strong solution. Indeed, since $c=27(9\times 2^{15/4})^4/16$ we get
$$
R_V = c^{-1/4}\nu {\lambda _1}^{1/2} =
c^{-\frac{1}{4}}=2/(9\times 2^{15/4}\times 27^{1/4})\approx 0.007.
$$
Nothing is known about regularity of solutions starting from initial
conditions lying in  a ball in $V$ with $R>R_V$, and so the
theoretical results currently available really are about ``small
data''. In this paper we present a method of numerical verification
of regularity for more general, physically reasonable initial
conditions. More precisely, denote by ${\mathscr B}_V(R)$ the ball
of radius $R$ in $V$.

We prove below under the assumption of regularity that  for any fixed value of $R>0$
the following statement can be verified
numerically in finite time.

\begin{statement}\label{st1} Every initial condition $u_0 \in {\mathscr B}_V(R)$ gives rise to a global strong
solution of problem (\ref{ns}).
\end{statement}

Of course, regularity
of solutions with initial conditions in some fixed ball in $V$ for
the non-dimensionalized form of the Navier-Stokes equations implies parallel results for other cubic domains and different values of the viscosity. Indeed, scaling initial condition changes its
enstrophy as we have
$$
\frac{\partial \tilde u}{\partial \tilde
x}=\frac{L^2}{4\pi ^2\nu} \frac{\partial u}{\partial x}.
$$
So
$$
\int_{[0,2\pi]^3} |D \tilde u(x)|^2\,\d\tilde x = \frac{L}{2\pi\nu
^2} \int_{[0,L]^3} |Du(x)|^2\,\d x
$$
and since $\lambda_1=\frac{4\pi
^2}{L^2}$ we see that
$$
||D \tilde u||^2 = \frac{1}{\nu^2
\lambda_1^{1/2}} ||Du||^2
$$
which is in agreement with theorem \ref{c}.

Thus a complete result on the regularity of solutions of the non-dimensionalized
Navier-Stokes equations would immediately imply the same result for
any $L$ and any $\nu$, while regularity only for $\|D\tilde u\|^2\le R$ in the non-dimensionalized
equation would provide regularity for all $\|Du\|^2\le2\pi R\nu^2/L$.

\section{Numerical verification of regularity for a single $u_0\in V$}\label{singleu0}

Before presenting a method of verifying regularity for all initial
conditions in bounded subsets of $H^2$ (and then $H^1$) we need a
result that allows us to check numerically the regularity of a
particular solution arising from some given initial condition. In
other words, we have to be able to make a judgement based only on
numerical calculations whether a given initial condition $u_0$ gives
rise to a regular solution $u$. The first step to this end is a
standard one: We construct the Galerkin approximations $u_n$ of
$u$ by solving the equation
$$
\frac{\d u_n}{\d
t}+Au_n+P_nB(u_n,u_n)=0
$$
with $u_n(0)=P_nu_0$, where $P_n$ is
orthogonal projection on the space spanned by first $n$
eigenfunctions of the Stokes operator.

We can deduce solely from the properties of these
Galerkin approximations whether they approximate a regular
solution on a finite time interval. A method for this was provided
by Chernyshenko et al.~(2006) under the assumption that $u_0\in H^m$ with $m\ge3$. Dashti \& Robinson (2006)
adapt these arguments to show that the same is possible given a regular
solution arising from an initial condition in $V$, and we now briefly
outline the main idea of this method.

The first ingredient is a robustness result for strong solutions:

\begin{theorem}\label{DRrobust}
There exists a constant $c>0$ such that if $u_0\in V$
gives rise to a strong solution $u(t)$ of (\ref{ournse}) on $[0,T^*]$ then so does any initial condition
$v_0$ and forcing $g(t)$ with
\begin{equation}\label{DR}
\|D(v_0-u_0)\|+\int_0^{T^*}\|Dg(s)\|\,\d s<c(T^*)^{-1/4}\exp\left(-c\int_0^{T^*}\|Du(s)\|^4+\|Du(s)\|\,\|Au(s)\|\,\d
s\right).
\end{equation}
(The constant $c$ is related to the constants arising in certain Sobolev embedding results.)
\end{theorem}

This has the following corollary, key to all that follows:

\begin{corollary}\label{DRcheck}
Suppose that $v\in L^\infty(0,T^*;V)\cap L^2(0,T^*;D(A))$ is a numerical approximation of $u$ such that
$$
\frac{\d v}{\d t}+Av+B(v,v)\in L^1(0,T^*;V)\cap L^2(0,T^*;H)
$$
and
\begin{eqnarray}
&&|Dv(0)-Du_0|+\int_0^{T^*}\left\|\frac{\d v}{\d t}(s)+Av(s)+B(v(s),v(s))\right\|_1\,\d s\nonumber\\
&&\qquad\qquad<c(T^*)^{-1/4}\exp\left(-c\int_0^{T^*}|Dv(s)|^4+|Dv(s)||Av(s)|\,\d
s\right).\label{givesreg}
\end{eqnarray}
Then $u$ is a regular solution of (\ref{ournse}) with $u\in L^\infty(0,T;V)\cap L^2(0,T;D(A))$.
\end{corollary}

Thus to verify that $u_0\in V$ gives rise to a strong solution on
$[0,T^*]$, it suffices to show that for sufficiently large $n$ the Galerkin
solution $u_n$ will satisfy (\ref{givesreg}) in Corollary
\ref{DRcheck}. If this is the case then on can compute ever more accurate
Galerkin approximations until (\ref{givesreg}) is eventually satisfied. If the approximated solution is regular then this method has
to be successful since it is the case that the Galerkin
approximations of regular solutions must converge in a sufficiently strong sense:

\begin{theorem}\label{DRgalerkin}
If $u_0\in V$ gives rise to a strong solution $u(t)$ on $[0,T^*]$
then the Galerkin approximations $u_n$ converge to $u$ strongly in
$L^\infty(0,T;V)\cap L^2(0,T;D(A))$.
\end{theorem}
\vspace*{0.15cm}

\noindent Of course, this neglects two issues: first, that the
numerical solution will not be a continuous function of $t$, but
given at discrete time points (cf.~comments in Chernyshenko et al.~(2006)). However, since Corollary
\ref{DRcheck} says nothing of the provenance of $v$, one can
construct $v$ via linear interpolation between these time points.
A further complication is that the integrals in (\ref{givesreg})
cannot be computed exactly and have to be approximated; but this
is not a major obstruction.

\section{Finite time of calculation}\label{finitettime}

Our goal in the following two sections is to prove, under the assumption of regularity, that the following statement may be verified numerically in a finite time.

\begin{statement}\label{stat1} Consider the non-dimensional 3d Navier-Stokes
equations on $[0,2\pi]^3$ with zero forcing,
\begin{equation}\label{problem} \frac{\partial u}{\partial t}-\nu\lap
u+(u\cdot\grad)u+\grad p=0\qquad\div u=0\qquad u(0)=u_0.
\end{equation} Then every initial condition in a ball in $V$ $$
{B}_{V}(R)=\{u_0:\ \|Du_0\|\le R\} $$ gives rise to a strong solution
that exists for all $t\ge0$. \end{statement}

Since the method in the previous section only works to verify the regularity of a given solution on some finite time interval, it is vital to show that this is all that is required.

\begin{lemma}\label{finitetime}
Any (weak) solution of (\ref{problem}) becomes regular after a
finite time $$ T^*=c^{1/2}\|u_0\|^2\le c^{1/2}\|Du_0\|^2, $$ where $c$
is the absolute constant from theorem \ref{c}.
\end{lemma}

\begin{proof}
Since $$ \int_0^T\|Du(s)\|^2\,\d s\le\|u_0\|^2, $$ it follows that on
any time interval $[0,T]$ with $T>T^*$ as above, there must exist
a $t_0$ such that $$ \|Du(t_0)\|^2\le c^{-1/2}, $$ and hence from
theorem \ref{c} it follows that the solution is regular for all $t\ge
T^*$.
\end{proof}

Thus if we are interested in regularity of a particular ball
${\mathscr B}_V(R)$ we have to check regularity of solutions only on some
fixed interval of time: now we can proceed with the proof of the
main theorem of this paper.

The most straightforward approach would be to verify the regularity
of solutions arising from some finite sequence of initial conditions
$u_{0n}$, where $1\le n \le N$ and $N$ is large enough. Then
Corollary \ref{DRcheck} would imply regularity of all initial
conditions in balls with radii $\delta_1, \delta_2,..., \delta_N $
centred at $u_{01}, u_{02},..., u_{0N}$ respectively. If these balls
cover all of ${\mathscr B}_R$ then the Statement \ref{stat1} would
be verified. However, several obstacles stand in the way of such a
happy ending. First of all, a ball in infinite dimensional space is
never compact and one cannot cover it with any finite number of
smaller balls. Even if we were interested only in checking
regularity in some finite-dimensional subspace of $V$  we would
still encounter a serious problem: without a proof of some general
uniform lower bound on the sequence $\delta _k$ we might not be able
to construct a proper covering (e.g.~it is not possible to cover the
unit interval with a sequence of intervals of length $\delta_n
=3^{-n}$).

To circumvent these difficulties we will prove in the next section an auxiliary result
under a stronger assumption on initial condition.

\section{Numerical verification of regularity for all initial conditions
in a ball in $H^2$}\label{H2ball}

In this section we show that it is possible to verify numerically
the following statement, which requires more regularity of the initial condition.

\begin{statement}\label{stat2} Consider the 3d Navier-Stokes
equations on $[0,2\pi]^3$ with zero forcing,
\begin{equation}\label{3dnse}
\frac{\partial u}{\partial t}-\lap
u+(u\cdot\grad)u+\grad p=0\qquad\div u=0\qquad u(0)=u_0.
\end{equation}
Then every initial condition in a ball $$ {\mathscr
B}_{H^2}(S)=\{u_0: |Au_0|\le S \}$$ gives rise to a strong
solution that exists for all $t\ge0$.\end{statement}

\begin{theorem}\label{VDA}
Statement \ref{stat2} can be verified numerically in a finite
time.
\end{theorem}

This theorem is proved via a series of subsidiary results. From
lemma \ref{finitetime} and the Poincar\'e inequality it follows
that in order to verify Statement \ref{stat2} we need only
prove regularity on the finite time interval $[0,T^*]$. It is
therefore very useful to have the following result, which gives
uniform bounds over the $H^1$ and $H^2$ norms of any solution with
$u_0\in {\mathscr
B}_{H^2}(S)$.

\begin{proposition}\label{FTregularity}
Suppose that Statement \ref{stat2} is true. Then for any $T>0$
there exist constants $D_S(T)$ and $E_S(T)$ such that for any
solution $u$ arising from an initial condition $u_0\in{\mathscr
B}_{H^2}(S)$ we have
\begin{equation}\label{bounds}
\sup_{0\le t\le T}\|Du(s)\|\le D_S\qqand
\int_0^{T}\|Au(s)\|^2\,\d s\le E_S.
\end{equation}
\end{proposition}
\begin{proof}
This result follows from the proof, but not the statement of Theorem
* in Foias \& Temam (19**), see also Theorem 12.10 in Robinson (2001). Here we
give a sketch of the relevant part of the argument. If the first
bound in (\ref{bounds}) does not hold then there must be sequences
$u_{0n}$ with $\|Au_{0n}\|\le S$ and $t_n\in[0,T]$ such that
$\|Du_n(t_n)\|\rightarrow\infty$ as $n\rightarrow\infty$. By
extracting appropriate subsequences we can assume that
$t_n\rightarrow T\in[0,T]$ and that $u_{0n}\rightharpoonup u_0$ in
$D(A)$: in particular, $\|Au_0\|\le S$ and $u_{0n}\rightarrow u_0$
in $H$. As in the standard proof of the existence of weak solutions,
one can show that the solutions $u_n(t)$ of
$$
\d u_n/\d t+Au_n+B(u_n,u_n)=0\qquad u(0)=u_{0n}
$$
are uniformly bounded in $L^\infty(0,T;H)\cap L^2(0,T;V)$; taking
the limit as $n\rightarrow0$ produces a weak solution $v$ of
$$
\d
v/\d t+Av+B(v,v)=0\qquad\mbox{with}\quad v(0)=v_0.
$$
However, by
assumption, this equation possesses a strong solution which is
unique in the class of weak solutions, and so in fact $v\in
L^\infty(0,T;V)\cap L^2(0,T;D(A))$. One can now use this
regularity of $v$ along with standard estimates on the equation
for the difference $v-u_n$ to show that $u_n\rightarrow v$
strongly in $L^2(0,T;V)$. Since $L^2$ convergence on an interval
implies that there exists a subsequence converging almost
everywhere, at each such point of convergence, $s$, we must have
$\|Du_n(s)\|\le S_0:=1+\|v\|_{L^\infty(0,T;V)}$ for $n$ sufficiently
large. Standard estimates can be used to show that there is a time
$\tau$ such that
$$
\|Du_n(s+t)\|\le 2(1+S_0)\qquad\mbox{for
all}\quad 0\le t\le\tau
$$
for all such $s$ and $n$. By covering
$[0,T]$ with a finite number of intervals $[s_j,s_j+\tau]$,
where $u_n(s_j)\rightarrow u(s_j)$ in $V$, one can deduce that in
fact $\|Du_n(t)\|\le 2(1+S_0)$ for all $t\in[0,T]$ and obtain a
contradiction. The second bound in (\ref{bounds}) then follows
from the inequality
$$
\frac{\d}{\d t}\|Du\|^2+\|Au\|^2\le c\|Du\|^6
$$
which can be obtained by standard methods.
\end{proof}

Now we are ready to prove that if an initial condition $u_0 \in
{\mathscr
B}_{H^2}(S)$ gives rise to a regular solution then there exist a
ball in $V$ of some radius $\delta$ depending on $S$ but not on
$u_0$ such that all initial conditions in that ball give rise to
regular solutions. In other words, we may say that if
Statement \ref{stat2} is true, then all points of ${\mathscr
B}_{H^2}(S)$ are
separated in $V$ from any initial conditions giving rise to a
solution with singularities (if such solutions exist) by at least some fixed
distance $\delta$.  Most importantly for us it means that there is a uniform lower bound
on the radius of  the ``balls of regularity'' whose existence is
guaranteed by theorem \ref{DRrobust}, centered at points of
${\mathscr
B}_{H^2}(S)$ .

\begin{proposition}\label{uniform}
Suppose that Statement \ref{stat2} is true. Then there exists a
$\delta=\delta(S)$ such that if $u_0 \in {\mathscr B}_{H^2}(S)$ then the value of $\delta$ required by the regularity check of Corollary \ref{DRcheck} is bounded below by some $\delta$: any initial condition  $v_0\in V$ with
$$
\|Du_0-Dv_0\|<\delta
$$
also gives rise to a strong solution on $[0,T^*]$.
\end{proposition}

\begin{proof}
  Given the result of Proposition \ref{FTregularity}, the integrand on the
right-hand side of (\ref{DR}) is bounded above:
$$
I:=\int_0^{T^*}\|Du(s)\|^4+\|Du(s)\|\,\|Au(s)\|\,\d s\le
I_S:=T^*D_S^4+(T^*D_S^2)^{1/2}E_S,
$$
and so the right-hand side
of (\ref{DR}) is bounded below,
$$
c\exp(-cI)\ge\delta:=c\exp(-cI_S).
$$
It follows that there exist a lower bound
on the radius of ``ball of regularity'' provided by theorem
\ref{DRrobust}
\end{proof}

Now we use a compactness argument. Loosely speaking, we
could cover ${\mathscr B}_{H^2}(S)$ by open balls of radius $\delta$ centred
at every point of ${\mathscr B}_{H^2}(S)$ and then choose finite sub-covering
(in topology of $V$). This argument, based on the Rellich-Kondrachov
theorem on the compact embedding of $H^2$ in $H^1$, is
formulated more precisely in the following lemma which also gives the
explicit form of points in ${\mathscr B}_{H^2}(S)$ that we need to test for regularity
numerically.

\begin{lemma}\label{cover}
Given $\delta>0$ there exist $N$, $M$ such that every $v_0\in
{\mathscr B}_{H^2}(S)$ can be approximated to within $\delta$ (in
the $V$ norm) by elements of the set $$ \mathcal{U}_{M,N} = \{ u_0:
u_0=\sum_{j=1}^N \alpha_jw_j,\hbox{ with }\alpha_j = a_j/2^M,\;
a_j\in {\mathbb Z};\ \|Au_0\|\le S\}$$ where $\{w_j\}$ are the
Stokes eigenfunctions.
\end{lemma}

\noindent Note that the proof shows that one can take $N$ such that
$\lambda_{N+1}\ge 2S^2/\delta$ and $M$ such that $2^{-M}<\delta/2$.

\begin{proof}
Take any $v \in {\mathscr B}_{H^2}(S)$. Since $\|Av\|\le S$, we have
\begin{eqnarray*}
\|D(P_n v -v)\|^2&=&\sum_{k=n+1}^\infty \lambda_k|(v,w_k)|^2\\
&\le& \lambda_{n+1}^{-1}\sum_{k=n+1}^\infty\lambda_k^2|(v,w_k)|^2\\\
 &\le&\lambda_{n+1}^{-1}\|Av\|^2\\
 &\le& \lambda_{n+1}^{-1}S^2.
 \end{eqnarray*}
Thus we can choose an $N$ that for any $n \geq N$ we have
$$
\|D(P_nv-v)\|^2 \le\lambda_{n+1}^{-1}S^2 \le
\frac{\delta}{2}
$$
and we can cover the finite dimensional ball $P_N ({\mathscr
B}_{H^2}(S))$ with a finite number of balls of radius $r \le
\frac{\delta}{2}$ with centres in $ \mathcal{U}_{M,N}$, where
$2^{-M}<\delta/2$. It follows that for any $v_0 \in {\mathscr
B}_{H^2}(S)$ we can find $\tilde u_0 \in \mathcal{U}_{M,N}$ such
that
$$
\|Dv_0-D\tilde u_0\| \le \|D(v_0-P_Nv_0)\|+ \|P_Nv_0-D\tilde u_0\|
\le \delta.
$$
\end{proof}

Finally we assemble all these ingredients to provide a proof of
theorem \ref{VDA}.

\begin{proof}[Proof of thereom \ref{VDA}]
Given $S>0$, let $\delta$ be the lower bound on the radius of the
``ball of regularity" provided by proposition \ref{uniform}. Then to
prove regularity for all initial conditions in ${\mathscr
B}_{H^2}(S)$ it suffices to check that the finite set of initial
conditions in ${\mathcal U}_{M,N}$ (the set defined in lemma
\ref{cover}) all give rise to strong solutions on $[0,T^*]$. This
can be done for each initial condition in turn using the procedure
outlined in Section \ref{singleu0}.
\end{proof}

Observe that the proof is based only on an assumption of regularity
of solutions arising from initial conditions in ${\mathscr
B}_{H^2}(S)$. Whether or not singularities may occur in solutions
arising from initial conditions $u_0$ with $\|Au_0\|>S$ has no
effect on the possibility of numerical verification of regularity in
the ball ${\mathscr B}_{H^2}(S)$.

\section{Numerical verification of regularity for all $u_0 \in {\mathscr B}_V(R)$}\label{H1ball}

Now we will use the results of previous section to prove the main
theorem of this paper

\begin{theorem}\label{mainheorem}
If the Navier-Stokes equations are regular, then statement
\ref{stat1} can be verified numerically in a finite time.
\end{theorem}

\begin{proof} We find a $T_R>0$ such that every solution
with initial condition $\|Du_0\|\le R$ is regular on the interval
$[0,T_R]$ and has an explicit bound on $\|Au(T_R)\|$, using a result
due to Foias \& Temam (1988) on the Gevrey regularity of solutions.

More explicitly, if we define
\begin{equation}
\e^{t A^{1/2}}= \sum_{n=0}^{\infty} \frac{t^n}{n!}A^{n/2}
\end{equation}
and set
$$
\tau=\frac{1}{K_1}(1+\|Du_0\|^2)^{-2},
$$
where $K_1\le 3266$ is an absolute constant depending on certain
algebraic and Sobolev embedding inequalities, Foias \& Temam (1988)
showed that for $0\le t\le 2\tau$, one has
$$
\|A^{1/2}\e^{tA^{1/2}}u(t)\|^2\le 2(1+\|Du_0\|^2).
$$

By expanding $A^{1/2}\e^{t A^{1/2}}u$ in terms of the eigenfunctions
of $A$ (cf.~proof of Lemma 2 in Friz \& Robinson, 2001) one can show
that
$$
\frac{(2\tau)^2}{2!}\|Au\|^2\le\|A^{1/2}\e^{\tau A^{1/2}}u\|,
$$
and hence that
$$
\|Au(\tau)\|^2\le K_1^2(1+\|Du_0\|^2)^5.
$$

So it suffices to check for regularity of initial conditions in some
ball of radius $S:=K_1(1+R^2)^{5/2}$ in $H^2$. Indeed, if any
initial condition $u_0\in {\mathscr B}_V(R)$ led to singularity so
would some initial condition $\tilde u_0 = u(\tau) \in {\mathscr
B}_{H^2}(S)$. So the regularity of solutions arising from ${\mathscr
B}_{H^2}(S)$ -- which  according to \ref{VDA} can be verified
numerically -- implies the regularity of solutions arising from
${\mathscr B}_V(R)$.\end{proof}

\vspace*{0.1cm}

Observe that in order to prove Statement
\ref{stat1} we need to assume that all solutions of the
Navier-Stokes equations are regular. Indeed, we cannot exclude the
possibility that some singularity may occur in a solution arising
from an initial condition $u_0$ with $\|Du_0\|>R$ but for which
$\|Au(\tau)\|< S$. Then, of course, we would not be able to
``verify'' numerically a false statement about regularity of
solutions arising from a ball ${\mathscr B}_{H^2}(S)$ in spite of the fact that
solutions arising from a ball ${\mathscr B}_{V}(R)$ would be regular.

\section{Conclusion}

We have shown that if the Navier-Stokes equations are regular on a cubic periodic domain, then regularity for any bounded set of initial conditions in $V$ (bounded enstrophy) can be verified numerically using an algorithm that is guaranteed to terminate. While this falls far short of a numerical proof of regularity, it is at least theoretically possible to check the validity of the equations for all initial conditions that one could hope to realise in experiment.

Of course, in this approach to regularity the outstanding open problem is to show that in fact regularity for all initial conditions in some bounded ball ${\mathscr B}_V(R)$ is in fact sufficient for regularity for all initial conditions, but this seems to be an extremely strong statement and it appears unlikely that this could be proved analytically without leading to a full proof of regularity.

\section*{Acknowledgments}
JCR is a Royal Society University Research Fellow, and would like
to thank the Society for all their support. WS is currently a
visiting fellow in the Mathematics Institute at the University of
Warwick under the Marie Curie Host Fellowship for the Transfer of
Knowledge.

\end{document}